\documentclass[11pt]{article}
\usepackage{amsfonts}
\usepackage{graphics}
\usepackage{indentfirst}
\usepackage{color}
\usepackage{cite}
\usepackage{latexsym}
\usepackage[paper=a4paper, left=2.1cm, right=2.1cm, top=2.2cm, bottom=1.5cm, headheight=5.5pt, footskip=0.8cm, footnotesep=0.8cm, centering, includefoot]{geometry}
\usepackage{amsmath}
\allowdisplaybreaks
\usepackage{amssymb}
\usepackage[dvips]{epsfig}
\usepackage{amscd}

\newtheorem{theorem}{Theorem}[section]
\newtheorem{remark}{Remark}[section]
\newtheorem{definition}{Definition}[section]
\newtheorem{lemma}{Lemma}[section]

\newtheorem{proposition}{Proposition}[section]

\DeclareMathOperator{\divv}{div}

\DeclareMathOperator{\trace}{tr}

\newcommand{\qn}{Q^{(n)}}
\newcommand{\qnm}{Q^{(n-1)}}
\newcommand{\pq}{L^p(0,t;L^q)}
\newcommand{\infq}{L^\infty(0,t;L^q)}
\newcommand{\infr}{L^\infty(0,t;L^r)}
\newcommand{\pinf}{L^p(0,t;L^\infty)}
\newcommand{\infinf}{L^\infty(0,t;L^\infty)}
\newcommand{\pr}{L^p(0,t;L^r)}

\makeatletter
\@addtoreset{equation}{section}
\makeatother
\makeatletter
\@addtoreset{equation}{section}
\makeatother

\title{Global solution to the three-dimensional liquid crystal flows of Q-tensor model}

\author{ Yao XIAO\thanks{Institute of Mathematical Sciences, The Chinese University of Hong Kong ({\tt yxiao@math.cuhk.edu.hk})
}
}

\date{ }

\begin{document}
	\maketitle
	
	\begin{abstract}
		We  study a complex non-Newtonian fluid that models the fluid of nematic liquid crystals. The system is a coupled system by a forced Navier-Stokes equation with a parabolic type equation of Q-tensors. In particular we invoke the maximal regularity of the Stokes equation and parabolic equation to show global strong solution for small initial data.
			\end{abstract}
	
	
	\section{Introduction}

Liquid Crystals present an intermediate state of matter between isotropic fluids and crystalline solids. It can be obtained typically by increasing the temperature of solid crystal (low-molecular-weight liquid crystals) or increase the concentration of some solvent (high-molecular-weight liquid crystals). In particular the nematic liquid crystals we are concerned with  in the current paper is composed of rod-like or disk-like molecules with head-to-tail symmetry and achiral property (identical to its mirror image). The molecules of nematic liquid crystals has short-range of positional order and long-range orientation order. More physical and chemical background can be found in the references in\cite{Gennes,Chan} and the references therein.

Several competing theories exist to capture this complexity of this non-Newtonian fluid. Ericksen and his collaborators initiated the study of the so-called Ericksen-Leslie-Prodi system \cite{Ericksen1,Ericksen-2,Leslie-1,Leslie-2} in 1960s, in which they used a unit vector to represent the macroscopic preferred orientation of molecules and the free energy is composed of kinetic energy and the Frank elastic energy. The corresponding static case is called Oseen-Frank theory \cite{HX}. For the simplified system we refer to \cite{Lin-1,Lin-2,Lin-3,Lin-4,Lin-5,Lin-6,Hong}. This theory is very successful in many aspects and matched with the experiments well especially during the phase intervals away from the threshold of phase transitions. However, the Ericksen-Leslie system could not predict the line defects commonly observed in the experiments. On the other hand, the unit director parameter could not describe the biaxial nematic liquid crystals, which means the alignment of the molecules has two preferred direction. To accommodate line defects and biaxial nematic liquid crystals, de Gennes proposed a traceless $d\times d$ symmetric matrix $Q$ as the new parameter, which is called $Q$-tensor. Suppose $x$ ia a point in the material domain, $\rho_x$ is a probability density function on molecular orientations which lie in $\mathbb{S}^{d-1}$ the unit sphere in $\mathbb{R}^d$. $\rho_x(\omega)=\rho_x(-\omega)$ due to the head-to-tail symmetry. The $Q$-tensor is the second-order moments of the probability measure on $\mathbb{S}^{d-1}$ with density $\rho_x$. 
\begin{equation*}
Q(x):=\int_{\mathbb{S}^{d-1}}(\omega\otimes\omega-\dfrac{I}{d})\rho_x(\omega)d \sigma(\omega).
\end{equation*}
It's easy to check that 
\begin{equation}
Q(x)\in S^d_{(0)}:=\{Q\in\mathbb{M}^{d\times d};Q=Q^T,\trace(Q)=0 \}.
\end{equation}

The term $d^{-1} I$ (which contains no information about the density $\rho_x$) is included by convention to render the $Q$-tensor identically zero when $\rho_x$ is the uniform distribution on the unit sphere. Thus $Q=0$ corresponds to the isotropic phase of the nematic liquid crystals. If $Q(x)$ admits two equal eigenvalues it corresponds uniaxial liquid crystals.  For the accordance of the $Q$-tensor system with the Ericksen-Leslie system, we refer to the work in \cite{MZ,ZZ}. If $Q(x)$ has three distinct eigenvalues, it's the biaxial liquid crystals.

Landau and de Gennes proposed the energy functional in terms of $Q$-tensor consists of the elastic energy and the bulk energy, one of the wildly accepted simplified form  which we adopt in the current paper is the following:
\begin{equation}\label{functional}
\mathcal{F}_{LG}:=\dfrac L2|\nabla Q|^2+\dfrac a2\trace(Q^2)-\dfrac b3\trace(Q^3)+\dfrac c4(\trace(Q^2))^2,
\end{equation} 
 where $L>0$, $a$, $b$ and $c$ are constants.
 
We also remark that concerning the well-posedness and the physicality of $Q$-tensor system, some progress on the modification of the energy functional has been made, see \cite{Ball-1,Ball-2,Mark}. 

For the dynamical model we use the one proposed by Beris and Edwards in \cite{Beris}.  The system writes as follows,

\begin{align}
\label{eqn1.1}
\begin{cases}
\partial_t Q+u\cdot \nabla Q-S(\nabla u,Q)=\Gamma H\\
\partial_t u+u\cdot \nabla u=\nu \Delta u+\nabla p+\divv (\tau+\sigma)\\
\divv u=0
\end{cases}
\end{align}
where $u(x,t)$ denotes velocity, $p$ is the pressure arising from the incompressibility. $\nu>0$ is the viscosity coefficient, and $\Gamma>0$ is the elasticity relaxation constant.  Moreover,
\begin{equation}
\label{def1.2}
S(\nabla u,Q):=(\xi D+\Omega)(Q+\dfrac{1}{d}Id)+(Q+\dfrac{1}{d}Id)(\xi D-\Omega)-2\xi(Q+\dfrac{1}{d}Id)tr(Q\nabla u)
\end{equation}
 in which $D:=\dfrac 12(\nabla u+\nabla^t u)$ and $\Omega:=\dfrac 12(\nabla u-\nabla^t u)$ are the symmetric and antisymmetric part, respectively, of the velocity gradient. $d$ is the dimension. $S(\nabla u, Q)$  describes how the flow gradient rotates and stretches the order-parameter $Q$. As proven by experiments and simple calculation in \cite{Beris}, a shear flow exert both tumbling and aligning effect on the liquid crystals, and $\xi$ denotes that ratio of these two effects, which depends on the environment and material parameters.

 We also denote:
\begin{equation}
\label{def1.3}
H:=-a Q+b[Q^2-\dfrac{Id}{d}tr(Q^2)]-cQtr(Q^2)+L\Delta Q,
\end{equation}

the symmetric part of the non-Newtonian (non-dissipative) stress:
\begin{align}
\label{def1.4}
\tau_{\alpha\beta}:=&-\xi(Q_{\alpha\gamma}+\dfrac{\delta_{\alpha\gamma}}{d})H_{\gamma\beta}-\xi H_{\alpha\gamma}(Q_{\gamma\beta}+\dfrac{\delta_{\gamma\beta}}{d})\\\nonumber
&+2\xi(Q_{\alpha\beta}+\dfrac{\delta_{\alpha\beta}}{d})Q_{\gamma\delta}H_{\gamma\delta}-L\partial_\beta Q_{\gamma\delta}\partial_{\alpha}Q_{\gamma\delta},
\end{align}
and the antisymmetric part:
\begin{equation}
\label{def1.5}
\sigma_{\alpha\beta}:=Q_{\alpha\gamma}H_{\gamma\beta}-H_{\alpha\gamma}Q_{\gamma\beta}.
\end{equation}

 Suppose $U$ is a bounded domain in $\mathbb{R}^3$ with smooth boundary in which the liquid crystals is confined, the following initial boundary condition is prescribed:
 \begin{equation}
 \label{1.8}
 (u,Q)|_{t=0}=(u_0,Q_0),\qquad (u,\partial_n Q)|_{\partial U}=(0,0),
 \end{equation}
 where $n$ is the outward normal vector on $U$.

As shown in\cite{PZ1}, when $\xi=0$, under the condition that $c>0$, the energy dissipation is shown through the particular cancelation of the system. Thus global weak solution in 2-D and 3-D follows, moreover in 2-D for general initial data, higher regularity and global strong solution is proved. In \cite{PZ2}, the authors generalize the conclusion to the cases $\xi$ sufficiently small under the same technique. The smallness condition can be removed if one restrict to a bounded domain.

Also in \cite{Ding}, they showed a global existence with a more general energy functional than $(\ref{functional})$. And global strong solution is proved for large viscosity too in 3-D .

In this paper we are interested in the global strong solution to (\ref{eqn1.1}) for suitable initial and boundary conditions in the 3-D case. By a strong solution we mean that (\ref{eqn1.1}) is satisfied almost everywhere. The global strong solution in 3-D fail in \cite{PZ1,PZ2} is mainly due to a key logarithm inequality fails to hold. Heuristically, the growth of higher-order norm is not well controlled in 3-D and the solution may break down in finite time. Roughly speaking, the linearization of $Q$ equation in (\ref{eqn1.1}) is a parabolic equation with polynomial source. In order to control the growth, we impose additional condition on the coefficients--which will be stated soon-- to let the source behave like a \textit{damping} term.

Recall the maximal regularity of Stokes equations in \cite{Danchin}, and the maximal regularity of parabolic equation in \cite{Amann},  also inspired by the work in \cite{HW} on the global existence of strong solution for simplified Ericksen-Leslie system with small initial data, we proved the global strong solution to (\ref{eqn1.1})  in 3-D for sufficiently small initial data. Moreover, we proved that this strong solution coincides with the weak solution constructed in \cite{PZ1,PZ2}.

The particular cancelation of higher order terms shown in \cite{PZ1} works well in the $L^2$ energy estimates. For the Cauchy problem, by the Littlewood-Paley decomposition, general $H^s$ estimates could also benefit from this cancelation. But now we want to invoke the maximal regularity of Stokes operator, we could only restrict to the bounded domain at the moment, as the Helmholtz decomposition may fail for general $L^p$ function when $p\neq 2$.  Without this cancelation, we have to shrink the solution space to find more regular solutions.

Another difficulty is that the polynomial term could result in exponential growth of the energy. This may not necessarily lead to the breakdown of the strong solution, but we still need to keep it under stronger control by the initial data so as to get uniform estimates for the whole time. Thus we make the following assumption:
\begin{equation}\label{damping}
ac\ge \dfrac 9{16}b^2,\quad a>0,\quad c>0.
\end{equation} 
so as to keep the lower order norm decay as shown in Proposition \ref{prop4.1}. Or one could just assume $a=0$, the result will be the same.

Moreover in this article we restrict ourselves to the case that $S(\nabla u, Q)=0$ due to technique reasons. From the modeling point of view that we consider those nematic liquid crystals that the stretching and rotating effects of the flow gradient can be neglected comparing to other dynamic behavior. Or equivalently we omit the terms with $\Delta Q$ in the  velocity equation, the difficulty and the method we adopt is basically the same. We are also uncomfortable about this assumption, but as it turns out our method could not control the gradient of velocity in the Q-tensor equation and $\Delta Q$ in the velocity equation at the same time. We hope to address this issue in the future work. 

Under this assumption, the entry form of system ({\ref{eqn1.1}}) reduces to the following:

\begin{align}
\label{eqn1.6}
\begin{cases}
 (\partial_t+u_\gamma\partial_\gamma)Q_{\alpha\beta}=\Gamma H_{\alpha\beta}\\
 (\partial_t+u_\beta\partial_\beta)u_\alpha=\nu\Delta u_\alpha+\partial_\alpha p-L\partial_{\beta}(\partial_\alpha Q_{\zeta\delta}\partial_\beta Q_{\zeta\delta})+L\partial_\beta(Q_{\alpha\gamma}\Delta Q_{\gamma\beta}-\Delta Q_{\alpha\gamma}Q_{\gamma\beta})\\
 \partial_\gamma u_\gamma =0
\end{cases}
\end{align}

Also, the quantity of the non-dimensional constants $\Gamma$, $\nu$ and $L$ do not play a role in our analysis, we will set them all to be $1$ for convenience.

The rest of the paper is organized as follows. In Sect. 2, our main results are stated, global existence of strong solution in 3-D for sufficiently small initial data. In Sect. 3, we recall the maximal regularity of parabolic equation and the Stokes equation and some useful embedding inequalities.  In Sect. 4, we proved the local existence of solution.  Finally in Sect. 5 we prove the global existence of solutions.
\newpage
\section{Main results}
We clarify the notations and conventions  first.

First of all a partial Einstein summation convention is adopted in the following, that is, we assume summation over repeated Greek indices $\alpha$, $\beta$, $\gamma$, $\dots $ but not the repeated Latin indices $i$, $j$, $k$, $\dots$, unless otherwise claimed.

In this article we use Frobenius norm of a matrix $|Q|:=\sqrt{\trace(Q^2)}=\sqrt{Q_{\alpha\beta}Q_{\alpha\beta}}$. And the Sobolev spaces, Besov spaces and the interpolation spaces of $Q$-tensors are all defined in terms of this norm. Let $|\nabla Q|^2=Q_{\alpha\beta,\gamma}Q_{\alpha\beta,\gamma}$ and $|\Delta Q|^2:=\Delta Q_{\alpha\beta}\Delta Q_{\alpha\beta}$, the higher order derivative norm is defined in the same pattern. We also denote $\nabla Q\odot\nabla Q$ a matrix in $\mathbb{R}^{d\times d}$ so that $(\nabla Q\odot\nabla Q)_{ij}:=\partial_i Q_{\alpha\beta}\partial_j Q_{\alpha\beta}$.

 To state the main results of this article, we define the appropriate function space where the solution is defined in.
 
 \begin{definition}
 	\label{def1.1}
 	For $T>0$ and $1<p,q,r<\infty$, $M^{p,q,r}_T$ denote the set of triplets $(u,Q,p)$ such that 
 	\begin{align*}
 	& u\in C([0,T],D^{1-\frac 1p,p}_{A_q})\cap L^p(0,T;W^{2,q}(U)\cap W^{1,q}_0(U)),\\
 	&\partial_t u\in L^p(0,T;L^q(U)),\qquad \divv u=0,\\
 	& Q\in C(0,T; B^{3-\frac 2p}_{r,p})\cap L^p(0,T;W^{3,r}(U)),\quad \partial_t Q\in L^p(0,T;W^{1,r}(U))\\
 	& p\in L^p(0,T;W^{1,q}(U)),\quad \int_{U} p dx=0.\\
 	& M^{p,q,r}_T \quad\text{is a Banach space with norm}\quad ||\cdot||_{M^{p,q,r}_T}.
 	\end{align*}
 \end{definition}
Notice that for the incompressible system, the pressure is still a solution up to a difference by constant, so the condition on $p$ in Definition \ref{def1.1} is satisfied automatically if we replace the pressure by
\begin{equation*}
p-\dfrac 1{|U|}\int_U pdx.
\end{equation*}
 
In the above definition, the space $D^{1-\frac 1p,p}_{A_q}$ is some fractional domain of the Stokes operator in $L^q$(\cite{Danchin}). Roughly speaking , the elements of this space are the vectors whose $2-\frac 2p$ derivatives are in $L^q$, divergence free, and vanish on $\partial U$. Moreover we have the following embedding
\begin{equation}
D^{1-\frac 1p,p}_{A_q}\hookrightarrow B^{2(1-\frac 1p)}_{q,p}\cap L^q.
\end{equation}
 While the Besov spaces on this bounded domain $U$ can be defined as the interpolation space as in \cite{Bergh}:
\begin{equation*}
B^{2(1-\frac 1p)}_{r,p}=(L^r,W^{1,r})_{1-\frac 1p,p}.
\end{equation*}
Now we state our existence result on the global strong solution. First we prove the local existence by an iterative method then obtain uniform estimates if the initial data is small enough.
 
 \begin{theorem}\label{existence}
 	Let $U$ be a bounded domain in $\mathbb{R}^3$ with $C^3$ boundary. Assume that $1< p<\frac 32$, $\frac 1p+\frac 1 q>1$, $r>q>3$, $q(1-\frac 1p)<\frac 32(\frac 1 q-\frac 1 r)$ and $u_0\in D^{1-\frac 1p,p}_{A_q}$, $Q_0\in B^{3-\frac 2p}_{r,p}\cap L^r$. Then
 	\begin{enumerate}
 		\item there exists a $T_0>0$, such that system (\ref{eqn1.6}) with initial-boundary conditions (\ref{1.8}) has a unique local strong solution $(u,Q,p)\in M^{p,q,r}_{T_0}$ in $U\times (0,T_0)$;
 		\item moreover, if (\ref{damping}) is satisfied, there exists a $\delta_0>0$, such that if the initial data satisfies
 		\begin{equation}\label{small}
 			||u_0||_{D_{A_q}^{1-\frac 1p,p}}\le \delta_0,\quad ||Q_0||_{B^{3-\frac 2p}_{r,p}}\le \delta_0,
 		\end{equation}
 		then (\ref{eqn1.6}),(\ref{1.8}) admits a unique global strong solution $(u,Q,p)\in M^{p,q,r}_T$ in $U\times(0,T)$ for all $T>0$.
 	\end{enumerate}
 \end{theorem}

\begin{remark}
	(i)The second part of the above theorem promises one a global solution near the equilibrium state $u=0$, $Q=0$,  which corresponds to the isotropic fluid state of the material. In a sense, we proved the stability of this equilibrium state of the liquid crystal at least for the model of (\ref{eqn1.6}).\\
	
	(ii) The set where $p$, $q$ and $r$ take values is not empty. For example,
	let $p=\frac {16}{15}$, $q=4$ and $r=15$, then they satisfy all the conditions. All these conditions ensure the following embedding:
	\begin{align}
	& D^{1-\frac 1p, p}_{A_q}\hookrightarrow B^{2(1-\frac 1p)}_{q,p}\hookrightarrow W^{2(1-\frac 1p),q}\hookrightarrow L^r;\label{2.3}\\
	& B^{3-\frac 2p}_{r,p}\hookrightarrow W^{3-\frac 2 p,r}\hookrightarrow W^{1,r}\hookrightarrow L^\infty.\label{2.4}
	\end{align}
\end{remark}

\vskip 0.5cm
\section{Preliminaries}

Let's recall the maximal regularities for the parabolic and Stokes operator, as well as some $L^\infty$ estimates. For the sake of statement, $C$ is the generic positive constant that may change from line to line.\\

Denote
\begin{equation*}
\mathcal{W}(0,T):= W^{1,p}(0,T;W^{1,r})\cap L^p(0,T;W^{3,r}).
\end{equation*}

First we state the maximal regularities for parabolic operator according to \cite{Amann} (p.p. 188, Remark 4.10.9(c)).
\begin{theorem}\label{thm2.1}
	Given $\omega_0\in B^{3-\frac 2p}_{r,p}$ and $f\in L^p(0,T;W^{1,r})$, the Cauchy problem
	\begin{equation*}
	\partial_t \omega-\Delta \omega= f,\qquad t\in(0,T),\qquad \omega(0)=\omega_0,
	\end{equation*}
	has a unique solution $\omega\in \mathcal{W}(0,T)$, and 
	\begin{equation*}
	||\omega||_{\mathcal{W}(0,T)}\le C(||f||_{L^p(0,T;W^{1,r})}+||\omega(0)||_{B^{3-\frac 2p}_{r,p}}),
	\end{equation*}
	where $C$ is independent of $\omega_0$, $f$ and $T$. Moreover, there exists a positive constant $c_0$ independent of $f$ and $T$ such that
	\begin{equation*}
	||\omega||_{\mathcal{W}(0,T)}\ge c_0\sup_{t\in(0,T)}||\omega(t)||_{B^{3-\frac 2p}_{r,p}}.
	\end{equation*}
\end{theorem}

%

Next we state the maximal regularities for the Stokes operator according to \cite{Danchin,Giga}.
\begin{theorem}\label{thm2.2}
	Let $U$ be a bounded domain with $C^{2+\epsilon}$ boundary in $\mathbb{R}^3$ and $1<p,q<\infty$. Assume that $u_0\in D^{1-\frac 1p,p}_{A_q}$ and $f\in L^p(0,\infty;L^q(U))$. Then the system
	\begin{align*}
	\begin{cases}
	\partial_t u-\Delta u+\nabla p= f,\\
	\int_{U} p=0,\quad \divv u=0,\\
	u|_{\partial U}=0,\quad u|_{t=0}=u_0.
	\end{cases}
	\end{align*}
	has a unique solution $(u,p)$ for all $T>0$ satisfying:
	\begin{equation}\label{3.1}
	||u(T)||_{D^{1-\frac 1p,p}_{A_q}}+||(\partial_t u,\Delta u,\nabla p)||_{L^p(0,T;L^q(U))}\le C(||u_0||_{D^{1-\frac 1p,p}_{A_q}}+||f||_{L^p(0,T;L^q(U))}).
	\end{equation}
	where $C=C(p,q,U)$.
\end{theorem}

\begin{remark}
	Notice that (\ref{3.1}) does not include the estimate for $||u||_{L^p(0,T;L^q)}$. However by the Dirichlet boundary condition $u|_{\partial U}=0$ Poincar$\acute{e}$ Inequality and Gagliardo-Nirenberg inequality imply:
	\begin{equation*}
	||u||_{L^{q}}\le  C||\nabla u||_{L^q}\le C||u||_{L^q}^{\frac 12}||\Delta u||_{L^q}^{\frac 12}\le \frac 12||u||_{L^q}+C||\Delta u||_{L^q}.
	\end{equation*}
\end{remark}
Thus(\ref{3.1}) can be rewritten as
\begin{equation*}
||u(T)||_{D^{1-\frac 1p,p}_{A_q}}+||(\partial_t u,u, \Delta u,\nabla p)||_{L^p(0,T;L^q(U))}\le C(||u_0||_{D^{1-\frac 1p,p}_{A_q}}+||f||_{L^p(0,T;L^q(U))}).
\end{equation*}

\begin{remark}
	While the similar conclusion for linear parabolic equations is true for both Cauchy problem and bounded domains, here for the Stokes equation, we only consider bounded domain with $C^{2+\epsilon}$ boundary. The main reason is that the Stokes operator is defined as the composition of Helmholtz projection operator (onto the solenoidal vector fields) and the Laplace Dirichlet operator (\cite{Giga2}), but for general unbounded domain, Helmholtz decomposition of $L^q$ functions may fail when $q\neq2$. Hence through out this paper we restrict ourselves to the bounded domain, and the corresponding boundary condition for $Q$ is mainly technical so that we can integrate by part without worrying the boundary terms.
\end{remark}

Next we introduce the $L^\infty$ estimate for the derivatives (cf. Lemma 4.1 \cite{Danchin}). Roughly speaking, the mismatch between space-time interpolation provides extra time decay. 
\begin{lemma}\label{lemma2.1}
	Let $1<p,q<\infty$ satisfy
	\begin{equation*}
	0<\dfrac p2 -\dfrac {3p}{2q}<1,
	\end{equation*}
	Then the following inequalities hold:
	\begin{align}
	||\nabla f||_{L^p(0,T;L^\infty)}\le CT^{\frac 12-\frac 3{2q}}||f||^{1-\theta}_{L^\infty(0,T:D^{1-\frac 1p,p}_{A_q})}||f||^\theta_{L^p(0,T;W^{2,q})},\label{3.2}
	\end{align}
	for some constant $C$ depending only on $U$, $p$, $q$ and 
	\begin{equation*}
	\dfrac{1-\theta}p=\dfrac 12-\dfrac 3{2q}.
	\end{equation*}
\end{lemma}

Similarly we also have
\begin{lemma}\label{lemma2.2}
 For {$r>3$, $0<\dfrac 1q-\dfrac 1r<\dfrac 23$ and $\dfrac 2p+\dfrac 3r>2$} we have  
	\begin{equation}\label{3.4}
	||\nabla^s f||_{L^p(0,T;L^\infty(U))}\le CT^{\frac{1-\theta}{p}}||f||_{L^\infty(0,T;B^{3-\frac 2p}_{r,p})}^{1-\theta}||f||^{\theta}_{L^p(0,T:W^{3,r}(U))}
	\end{equation}
	with 
	\begin{equation*}
	\dfrac{1-\theta}{p}=\dfrac 1{2}(3-s-\dfrac 3r), \quad\text{where}\quad \dfrac p{2}(3-s-\dfrac 3 r)\in (0,1), \quad s=1, 2.
	\end{equation*}
	\begin{equation}{\label{n3.5}}
	||\nabla^3 f||_{L^p(0,T; L^q(U))}\le CT^{\frac {1-\theta'}{p}}||f||_{L^\infty(0,T;W^{1,r})}^{1-\theta'}||f||_{L^p(0,T;W^{3,r})}^{\theta'},
	\end{equation}
	with
	\begin{equation*}
	\dfrac{1-\theta'}{p}=\dfrac 3{2p}(\dfrac 1q-\dfrac1r),\quad\text{where}\quad \dfrac 32(\dfrac 1q-\dfrac 1r)\in (0,1).
	\end{equation*}
\end{lemma}

Notice that under the conditions of $p$, $q$ and $r$ in Theorem \ref{existence}, the requirement here is naturally satisfied.

\noindent
\textit{Proof}:\\

First notice the embedding $B^0_{\infty,1}\hookrightarrow L^\infty$ and interpolation relation
\begin{equation*}
(B^{3-s-\frac 2p-\frac 3r}_{\infty,\infty}, B^{3-s-\frac 3r}_{\infty,\infty})_{\theta,1}=B^0_{\infty,1}\quad\text{with}\quad\dfrac{1-\theta}{p}=\dfrac 12(3-s-\dfrac 3r)
\end{equation*}

Thus one have
\begin{equation*}
||\nabla^s f||_{L^\infty}\le C||\nabla^s f||_{B^0_{\infty,1}}\le C||\nabla^s f||_{B^{3-s-\frac 2p-\frac 3r}_{\infty,\infty}}^{1-\theta}||\nabla^s f||_{B^{3-s-\frac 3r}_{\infty,\infty}}^\theta.
\end{equation*}

Then invoke the following embedding relation
\begin{equation*}
B^{3-\frac 2p}_{r,p}\hookrightarrow B^{3-\frac 2p-\frac 3r}_{\infty,\infty},\quad W^{3,r}\hookrightarrow B^{3}_{r,\infty}\hookrightarrow B^{3-\frac 3r}_{\infty,\infty}
\end{equation*}
one obtain
\begin{align*}
||\nabla^s f||_{L^p(0,T;L^\infty)}\le & C(\int_0^T||\nabla^s f||_{B^{3-s-\frac 2p-\frac 3r}_{\infty,\infty}}^{p(1-\theta)}||\nabla^s f||_{B^{3-s-\frac 3r}_{\infty,\infty}}^{p\theta}dt)^{\frac 1p}\\
\le& C (\int_0^T||f||_{B^{3-\frac 2p}_{r,p}}^{p(1-\theta)}||f||_{W^{3,r}}^{p\theta}dt)^{\frac 1p}\\
\le& CT^{\frac{1-\theta}p}||f||_{L^\infty(0,T;B^{3-\frac 2p}_{r,p})}^{1-\theta}||f||_{L^p(0,T;W^{3,r})}^\theta
\end{align*}
Thus \eqref{3.4} is proved.

By Gagliardo-Nirenberg interpolation inequality:
\begin{equation*}
||\nabla^3 f||_{L^q}\le C||\nabla f||_{L^r}^{1-\theta'}||\nabla^3 f||_{L^r}^{\theta'}
\end{equation*}
where $1-\theta'=\dfrac 32(\dfrac 1q-\dfrac 1 r)$. Using H$\ddot{o}$lder inequality as above one get \eqref{n3.5}
\hfill
$\Box$

\begin{remark}
\eqref{n3.5} is mainly used to control $\Delta Q$ appeared in the velocity equation, this is possible as we choose $r>q$, but the sacrifice is that we can no longer control $\nabla u$ in the Q-tensor equation. If one choose $r<q$ it's the other way around, we could control $\nabla u$ while $\Delta Q$ is lost. What puzzles us is how can we control these two at the same time.
\end{remark}

\vskip 0.5 cm
\section{Local existence of solution}

In this section, we prove the existence  and the uniqueness of local strong solution by an iterative method, which is divided into several steps, including constructing the approximation scheme, obtaining the uniform estimate, showing the convergence, consistency and uniqueness of the limit.\\

\noindent
\textit{Step 1. Construction of the approximate solution}.\\

We initialize the construction of the approximation by setting $Q^{(0)}:=Q_0$ and $u^0:=u_0$. For given $(Q^{(n)},u^n,p^n)$ where $n\in \mathbb N$, define $(Q^{(n+1)},u^{n+1},p^{n+1})$ as the global solution of the following approximation system:
\begin{align}\label{eqn3.1}
\begin{cases}
&\partial_t u^{n+1}-\Delta u^{n+1}+\nabla p^{n+1}=-u^n\cdot\nabla u^n-\divv(\nabla Q^{(n)}\odot\nabla Q^{(n)}+\Delta Q^{(n)}Q^{(n)}-Q^{(n)}\Delta Q^{(n)})\\
&\partial_t Q^{(n+1)}-\Delta Q^{(n+1)}=-u^n\cdot\nabla Q^{(n)}-aQ^{(n)}+b[{Q^{(n)}}^2-\dfrac{Id}{d}tr({Q^{(n)}}^2)]-cQ^{(n)}tr({Q^{(n)}}^2)\\
&\divv u^{n+1}=0,\qquad \int_{U}p^{n+1}dx=0.\\
&(u^{n+1},Q^{(n+1)})|_{t=0}=(u_0,Q_0),\qquad (u^{n+1},\partial_nQ^{(n+1)})|_{\partial U}=(0,0).
\end{cases}
\end{align} 

Apparently such existence and uniqueness of the solution to the linear Stokes equation and linear parabolic equation ensures that $(u^{n+1},Q^{(n+1)},p^{n+1})$ is well-defined. Moreover by Theorem \ref{thm2.1} and \ref{thm2.2}, there exists a sequence of $\{(u^n,Q^{(n)},p^n)\}_{n\in\mathbb{N}}\in M^{p,q,r}_T$ for all positive $T$.\\

\noindent
\textit{Step 2. Uniform estimate for small fixed time $T$}.\\

The purpose of this step is to find certain positive time $T$ independent of $n$.

According Theorem \ref{thm2.1}, the maximal regularities of parabolic operator,
\begin{equation}\label{eqn3.2}
\begin{aligned}
&||Q^{(n+1)}||_{L^\infty(0,T;B^{3-\frac 2p}_{r,p})}+||Q^{(n+1)}||_{\mathcal{W}(0,T)}\lesssim ||Q_0||_{B^{3-\frac 2p}_{r,p}}+\underbrace{||u^n\cdot\nabla Q^{(n)}||_{L^p(0,T;W^{1,r}(U))}}_{\mathcal{I}_1}\\
&+\underbrace{||aQ^{(n)}+b[{Q^{(n)}}^2-\dfrac {Id}dtr({Q^{(n)}}^2)]-cQ^{(n)}tr({Q^{(n)}}^2)||_{L^p(0,T;W^{1,r}(U))}}_{\mathcal{I}_2}
\end{aligned}
\end{equation}

Similarly by Theorem \ref{thm2.2}, the maximal regularities of Stokes operator,
\begin{equation}\label{eqn3.3}
\begin{aligned}
&||u^{n+1}(T)||_{D^{1-\frac 1p,p}_{A_q}}+||(\nabla p^{n+1},u^{n+1},\Delta u^{n+1},\partial_t u^{n+1})||_{L^p(0,T;L^q(U))}\lesssim ||u_0||_{D^{1-\frac 1p,p}_{A_q}}\\
&+\underbrace{||u^n\cdot\nabla u^n||_{L^p(0,T;L^q(U))}}_{\mathcal{J}_1}+\underbrace{||\nabla\cdot(\nabla Q^{(n)}\odot\nabla Q^{(n)})||_{L^p(0,T;L^q(U))}}_{\mathcal{J}_2}\\
&+\underbrace{||\nabla\cdot(Q^{(n)}\Delta Q^{(n)}-\Delta Q^{(n)}Q^{(n)})||_{L^p(0,T;L^q(U))}}_{\mathcal{J}_3}
\end{aligned}
\end{equation}

Now define
\begin{align*}
U^n(t)=&||u^n||_{L^\infty(0,t;D^{1-\frac 1p,p}_{A_q})}+||u^n||_{L^p(0,t;W^{2,q})}+||\partial_t u^n||_{L^p(0,t;L^q)}\\
&+||\nabla p||_{\pq}+||Q^{(n)}||_{L^\infty(0,t:B^{3-\frac 2p}_{r,p})}+||Q^{(n)}||_{\mathcal{W}(0,t)}.
\end{align*}

With the help of Lemma \ref{lemma2.1} and Lemma \ref{lemma2.2}, we then control the right hand side of (\ref{eqn3.2}) and (\ref{eqn3.3}) by $U^n(t)$, obtaining the inequality of $U^n(t)$. And in the following estimates $\theta$ may denote different values in different terms, and we only compute the leading order terms as the lower order terms will trivially follow.

\begin{align*}
\mathcal{I}_1\lesssim &||u^n||_{L^\infty(0,t;L^r)}||\nabla^2 Q^{(n)}||_{L^p(0,t;L^\infty)}+||\nabla u||_{L^p(0,t;L^\infty)}||\nabla Q^{(n)}||_{L^\infty(0,t;L^r)}\\
\lesssim & T^{\frac 12(1-\frac 3r)}||u^n||_{L^\infty(0,t;L^r)}||Q^{(n)}||^{1-\theta}_{L^\infty(0,t;B^{3-\frac 2p}_{r,p})}||Q^{(n)}||_{L^p(0,t;W^{3,r})}^\theta\\
& +T^{\frac 12(1-\frac 3q)}||u||_{L^\infty(0,t;D^{1-\frac 1p,p}_{A_q})}^{1-\theta}||u||_{L^p(0,t;W^{2,q})}||Q^{(n)}||_{L^\infty(0,t;W^{1,r})}\\
\lesssim & (T^{\frac 12(1-\frac 3r)}+T^{\frac 12(1-\frac 3q)})[U^n(t)]^2,
\end{align*}
where we have used \eqref{3.2}, (\ref{3.4}) for $s=2$, and \eqref{2.3}\eqref{2.4}.


\begin{align*}
\mathcal{I}_2\lesssim & |a|||Q^{(n)}||_{L^p(0,t;W^{1,r})}+|b|||[Q^{(n)}]^2||_{L^p(0,t;W^{1,r})}+c||[Q^{(n)}]^3||_{L^p(0,t;W^{1,r})}\\
\lesssim & T^\frac 1p||Q^{(n)}||_{L^\infty(0,t;L^r)}(1+||Q^{(n)}||_{L^\infty(0,t;L^\infty)}+||Q^{(n)}||_{L^\infty(0,t;L^\infty)}^2)\\
& +||\nabla Q^{(n)}||_{L^p(0,t;L^\infty)}(1+||Q^{(n)}||_{L^\infty(0,t;L^r)}+||Q^{(n)}||_{L^\infty(0,t;L^r)}||Q^{(n)}||_{L^\infty(0,t;L^\infty)})\\
\lesssim & (T^\frac 1p+T^{\frac 12(2-\frac 3r)} )U^n(t)(1+U^n(t)+[U^n(t)]^2).
\end{align*}
where we have used \eqref{2.4}

Similarly for the right hand side of (\ref{eqn3.3}),

\begin{align*}
\mathcal{J}_1\lesssim &||u^n||_{L^\infty(0,t;L^q)}||\nabla u^n||_{L^p(0,t;L^\infty)}\\
\lesssim & T^{\frac 12(1-\frac 3q)}||u||_{L^\infty(0,t;L^q)}||u^n||_{L^\infty(0,t;D^{1-\frac 1p,p}_{A_q})}^{1-\theta}||u^n||_{L^p(0,t;W^{2,q})}^\theta\\
\lesssim & T^{\frac 12(1-\frac 3q)}[U^n(t)]^2,
\end{align*}
where we have used $D^{1-\frac 1p,p}_{A_q}\hookrightarrow L^q$ and (\ref{3.2}).

\begin{align*}
\mathcal{J}_2\lesssim &2||\nabla^2 Q^{(n)}||_{L^p(0,t;L^\infty)}||\nabla Q^{(n)}||_{L^\infty(0,t;L^q)}\\
\lesssim & T^{\frac 12(1-\frac 3 r)}||Q^{(n)}||_{L^\infty(0,t;B^{3-\frac 2p}_{r,p})}^{1-\theta}||Q^{(n)}||_{L^p(0,t;W^{3,r})}^\theta||Q^{(n)}||_{L^\infty(0,t;B^{3-\frac 2p}_{r,p})}\\
\lesssim & T^{\frac 12(1-\frac 3r)}[U^n(t)]^2,
\end{align*}
where we have used (\ref{3.4}) for $s=2$ and \eqref{2.4}.

The following is the highest order term, and the key in our estimates,
\begin{align*}
\mathcal{J}_3\lesssim &||\nabla Q^{(n)}\Delta Q^{(n)}||_{L^p(0,t;L^q)}+||Q^{(n)}\nabla\Delta Q^{(n)}||_{L^p(0,t;L^q)}\\
\lesssim & ||\nabla Q^{(n)}||_{L^\infty(0,t;L^q)}||\Delta Q^{(n)}||_{L^p(0,t;L^\infty)}+||\nabla\Delta Q^{(n)}||_{L^p(0,t;L^q)}||Q^{(n)}||_{L^\infty(0,t;L^\infty)}\\
\lesssim & T^{\frac 12(1-\frac 3r)}||Q^{(n)}||_{L^\infty(0,t;B^{3-\frac 2p}_{r,p})}^{1-\theta}||Q^{(n)}||_{L^p(0,t;W^{3,r})}^\theta||Q^{(n)}||_{L^\infty(0,t;W^{1,q})}\\
&+T^{\frac 3{2p}(\frac 1q-\frac 1r)}||Q^{(n)}||_{L^\infty(0,t;W^{1,r})}^{1-\theta'}||Q^{(n)}||_{L^p(0,t;W^{3,r})}^{\theta'}||Q^{(n)}||_{L^\infty(0,t;L^\infty)}\\
\lesssim & (T^{\frac 12(1-\frac 3r)}+T^{\frac 3{2p}(\frac 1q-\frac 1r)})[U^n(t)]^2.
\end{align*}
where we have used (\ref{3.4}) for $s=2,3$, \eqref{n3.5} and \eqref{2.4}.

Adding up (\ref{eqn3.2}) and (\ref{eqn3.3}), substituting all the estimates above, we have for all $t\in [0,T]$,
\begin{equation}\label{eqn3.4}
\begin{aligned}
U^{n+1}(t)\le & C(T^{\frac 12(1-\frac 3 r)}+T^{\frac 12(1-\frac 3q)}+T^{\frac 3{2p}(\frac 1 q-\frac 1 r)})[U^n(t)]^2\\
& +C(T^\frac 1pU^n(t)+T^{\frac 12(2-\frac 3 r)}[U^n(t)]^3)+U_0
\end{aligned}
\end{equation}

Suppose $U^n(t)\le 6C U_0$ on $[0,T]$ with 
\begin{equation}\label{eqn3.5}
\begin{aligned}
& 0<T\le (\dfrac 5{6C+108C^2U_0+216C^3U_0^2})^{\frac{1}{\alpha}}\le 1,\quad \text{or}\\
& 1<T\le (\dfrac 5{6C+108C^2U_0+216C^3U_0^2})^{\frac 1\beta}.
\end{aligned}
\end{equation}
where $\alpha=\min\{\frac 12(1-\frac 3 q),\frac 1{2p}(\frac 1q-\frac 1r),\frac 1p,\frac 12(2-\frac 3r)\}$, and $\beta=\max\{\frac 12(1-\frac 3 q),\frac 1{2p}(\frac 1q-\frac 1r),\frac 1p,\frac 12(2-\frac 3r)\}$.

Then direct computation implies
\begin{equation*}
U^{n+1}(t)\le 6 CU_0 \quad \text{on}\quad [0,T].
\end{equation*}

Relax unitil $C>\frac 16$, the assumption is satisfied for $n=1$, then by induction we have proved the following lemma,
\begin{lemma}\label{lemma3.1}
	For all $t\in [0,T_0]$ with $T_0$ satisfying (\ref{eqn3.5}),
	\begin{equation}\label{eqn3.6}
	U^n(t)\le 6CU_0, \quad \forall n\in \mathbb N.
	\end{equation}
\end{lemma}

\noindent
\textit{Step 3. Convergence of the approximation solution}.\\

$M^{p,q,r}_T$ is a Banach space, we now prove
\begin{lemma}\label{lemma3.2}
	$\{(u^n,Q^{(n)},p^n)\}$ is a Cauchy sequence and converges in $M^{p,q,r}_{T_0}$.
\end{lemma}
\textit{Proof}: Let
\begin{equation*}
\delta u^n:=u^{n+1}-u^n,\quad \delta Q^{(n)}:=Q^{(n+1)}-Q^{(n)},\quad \delta p^n:=p^{n+1}-p^n.
\end{equation*}
Define
\begin{equation*}
\begin{aligned}
\delta U^n(t):=&||\delta u^n||_{L^\infty(0,t;D^{1-\frac 1p,p}_{A_q})}+||\delta u^n||_{L^p(0,t;W^{2,q})}+||\partial_t\delta u^n||_{L^p(0,t;L^q)}\\
&+||\nabla\delta p||_{\pq}+||\delta Q^{(n)}||_{L^\infty(0,t;B^{3-\frac 2p}_{r,p})}+||\delta Q^{(n)}||_{\mathcal{W}(0,t)}.
\end{aligned}
\end{equation*}
Then the triplet $(\delta u^n,\delta Q^{(n)},\delta p^n)$ satisfies 
\begin{align}\label{eqn3.7}
\begin{cases}
&\partial_t\delta u^n-\Delta \delta u^n+\nabla \delta p^n=-u^n\cdot\nabla u^n+u^{n-1}\cdot\nabla u^{n-1}-\divv(\nabla Q^{(n)}\odot\nabla Q^{(n)}-\nabla Q^{(n-1)}\odot\nabla Q^{(n-1)})\\
&-\divv(\Delta Q^{(n)}Q^{(n)}-Q^{(n)}\Delta Q^{(n)}-\Delta Q^{(n-1)}Q^{(n-1)}+Q^{(n-1)}\Delta Q^{(n-1)})\\
&\partial_t\delta Q^{(n)}-\Delta\delta Q^{(n)}=-u^n\cdot\nabla Q^{(n)}+u^{n-1}\cdot\nabla Q^{(n-1)}-a\delta Q^{(n-1)}\\
&+b[{Q^{(n)}}^2-{Q^{(n-1)}}^2-\dfrac{Id}{d}(\trace({Q^{(n)}}^2)-\trace({Q^{(n-1)}}^2))]-cQ^{(n)}\trace({Q^{(n)}}^2)+cQ^{(n-1)}\trace({Q^{(n-1)}}^2),\\
&\divv \delta u^n=0,\qquad \int_{U}\delta p^ndx=0,\\
&(\delta u^n,\delta Q^{(n)})|_{t=0}=(0,0),\qquad (\delta u^n,\partial_n\delta Q^{(n)})|_{\partial U}=(0,0).
\end{cases}
\end{align}
By the same token of step 2, by using Lemma \ref{lemma2.1}, \ref{lemma2.2} and Lemma \ref{lemma3.1}.
\begin{align}
&||-u^n\cdot\nabla u^n+u^{n-1}\cdot\nabla u^{n-1}||_{L^p(0,t;L^q)}\nonumber\\
=&||-\delta u^{n-1}\cdot\nabla u^n-u^{n-1}\cdot\nabla\delta u^{n-1}||_{L^p(0,t;L^q)}\nonumber\\
\le& ||\delta u^{n-1}||_{L^\infty(0,t;L^q)}||\nabla u^n||_{L^p(0,t;L^\infty)}+||u^{n-1}||_{L^\infty(0,t;L^q)}||\nabla\delta u^{n-1}||_{L^p(0,t;L^\infty)}\nonumber\\
\lesssim & T^{\frac 12-\frac 3{2q}}U^n(t)||\delta u^{n-1}||_{L^\infty(0,t;L^q)}+U^{n-1}(t)||\nabla\delta u^{n-1}||_{L^p(0,t;L^\infty)}\nonumber\\
\lesssim& 6CU_0(T^{\frac 12-\frac 3{2q}}||\delta u^{n-1}||_{L^\infty(0,t;L^q)}+||\nabla\delta u^{n-1}||_{L^p(0,t;L^\infty)})\nonumber\\
\lesssim &T^{\frac 12-\frac 3{2q}}\delta U^{n-1}(t),\label{3.8}\\
\nonumber\\
&||-\divv(\nabla Q^{(n)}\odot\nabla Q^{(n)}-\nabla Q^{(n-1)}\odot\nabla Q^{(n-1)})||_{L^p(0,t;L^q)}\nonumber\\
=&||\divv (\nabla\delta Q^{(n-1)}\odot\nabla Q^{(n)}+\nabla Q^{(n-1)}\odot\nabla\delta Q^{(n-1)})||_{L^p(0,t;L^q)}\nonumber\\
 \le&||\divv(\nabla\delta Q^{(n-1)}\odot\nabla Q^{(n)})||_{L^p(0,t;L^q)}+||\divv(\nabla Q^{(n-1)}\odot\nabla\delta Q^{(n-1)})||_{L^p(0,t;L^q)}\nonumber\\
\lesssim & ||\nabla^2\delta Q^{(n-1)}||_{L^p(0,t;L^\infty)}||\nabla Q^{(n)}||_{L^\infty(0,t;L^q)}+||\nabla\delta Q^{(n-1)}||_{L^\infty(0,t;L^q)}||\nabla^2 Q^{(n)}||_{L^p(0,t;L^\infty)}\nonumber\\
&+||\nabla^2 Q^{(n-1)}||_{L^p(0,t;L^\infty)}||\nabla\delta Q^{(n-1)}||_{L^\infty(0,t;L^q)}+||\nabla Q^{(n-1)}||_{L^\infty(0,t;L^q)}||\nabla^2\delta Q^{(n-1)}||_{L^p(0,t;L^\infty)}\nonumber\\
\lesssim & (||\nabla^2\delta Q^{(n-1)}||_{L^p(0,t;L^\infty)}+||\nabla\delta Q^{(n-1)}||_{L^\infty(0,t;L^q)}T^{\frac 12(1-\frac 3r)})(U^n(t)+U^{n-1}(t))\nonumber\\
\lesssim& 6CU_0(||\nabla^2\delta Q^{(n-1)}||_{L^p(0,t;L^\infty)}+||\nabla\delta Q^{(n-1)}||_{L^\infty(0,t;L^q)}T^{\frac 12(1-\frac 3r)})\nonumber\\
\lesssim & T^{\frac 12(1-\frac 3 r)}\delta U^{n-1}(t),\label{3.9}\\
\nonumber\\
&||-\divv(\Delta Q^{(n)}Q^{(n)}-\qn\Delta\qn-\Delta\qnm\qnm+\qnm\Delta\qnm)||_{L^p(0,t;L^q)}\nonumber\\
\le&||\divv(\Delta\delta\qnm\qn+\Delta\qnm\delta\qnm)||_{\pq}\nonumber\\
&+||\divv(\delta\qnm\Delta\qn+\qnm\Delta\delta\qnm)||_{\pq}\nonumber\\
\le&||\nabla\Delta\delta\qnm||_{\pq}(||\qn||_{\infinf}+||\qnm||_{\infinf})\nonumber\\
&+||\Delta\delta\qnm||_{\pinf}(||\nabla\qn||_{\infq}+||\nabla\qnm||_{\infq})\nonumber\\
&+||\delta\qnm||_{\infinf}(||\nabla\Delta\qn||_{\pq}+||\nabla\Delta\qnm||_{\pq})\nonumber\\
&+||\nabla\delta\qnm||_{\infq}(||\Delta\qnm||_{\pinf}+||\Delta\qn||_{\pinf})\nonumber\\
\lesssim & (U^n(t)+U^{n-1}(t))(||\nabla\Delta\delta\qnm||_{\pq}+||\Delta\delta\qnm||_{\pinf}\nonumber\\
&+T^{\frac12(1-\frac 3r)}||\nabla\delta\qnm||_{\infq}+T^{\frac 3{2p}(\frac 1q-\frac 1r)}||\delta\qnm||_{\infinf})\nonumber\\
\lesssim& 6CU_0(||\nabla\Delta\delta\qnm||_{\pq}+||\Delta\delta\qnm||_{\pinf}\nonumber\\
&+T^{\frac 12(1-\frac 3r)}||\nabla\delta\qnm||_{\infq}+T^{\frac 3{2p}(\frac 1q-\frac 1r)}||\delta\qnm||_{\infinf})\nonumber\\
\lesssim &(T^{\frac 3{2p}(\frac 1q-\frac 1r)}+T^{\frac 12(1-\frac 3r)})\delta U^{n-1}(t),\label{3.10}\\
\nonumber\\
&||-u^n\cdot\nabla\qn+u^{n-1}\cdot\nabla\qnm||_{L^p(0,t;W^{1,r})}   \nonumber \\
\lesssim & ||\delta u^{n-1}||_{\infr}||\nabla\qn||_{\pinf}+||u^{n-1}||_{\infr}||\nabla\delta\qnm||_{\pinf}    \nonumber\\
&+||\nabla\delta u^{n-1}||_{\pinf}||\nabla \qn||_{\infr}+||\delta u^{n-1}||_{\infr}||\nabla^2\qn||_{\pinf} \nonumber\\
&+||\nabla u^{n-1}||_{\pinf}||\nabla\delta\qnm||_{\infr}+||u^{n-1}||_{\infr}||\nabla^2\delta\qnm||_{\pinf}  \nonumber\\
\lesssim & 6CU_0\large[(T^{\frac12(2-\frac 3r)}+T^{\frac 12(1-\frac 3r)})||\delta u^{n-1}||_{\infr}+||\nabla\delta u^{n-1}||_{\pinf} \nonumber\\
&+||\nabla\delta\qnm||_{\pinf}+T^{\frac 12(1-\frac 3q)}||\nabla\delta\qnm||_{\infr}+||\nabla^2\delta\qnm||_{\pinf}\large] \nonumber\\
\lesssim& (T^{\frac 12(2-\frac 3r)}+T^{\frac 12(1-\frac 3r)}+T^{\frac 12(1-\frac 3q)})\delta U^{n-1}(t),\label{3.11}\\
&\nonumber\\
&|b|||{\qn}^2-{\qnm}^2-\dfrac{Id}{d}(\trace({\qn}^2)-\trace({\qnm}^2)||_{L^p(0,t;W^{1,r})}\nonumber\\
\lesssim&||\delta\qnm||_{\infr}(||\qn||_{\pinf}+||\qnm||_{\pinf})\nonumber\\
&+||\nabla\delta\qnm||_{\pinf}(||\qn||_{\infr}+||\qnm||_{\infr})\nonumber\\
&+||\delta\qnm||_{\infr}(||\nabla\qn||_{\pinf}+||\nabla\qnm||_{\pinf}) \nonumber\\
\lesssim &6CU_0((T^{\frac 1p}+T^{\frac 12(2-\frac 3r)})||\delta\qnm||_{\infr}+||\nabla\delta\qnm||_{\pinf})\nonumber\\
\lesssim & (T^\frac 1p+T^{\frac 12(2-\frac 3r)})\delta U^{n-1}(t)\label{3.13}\\
 \nonumber\\
&|c|||\qn\trace({\qn}^2)-\qnm\trace({\qnm}^2)||_{L^p(0,t;W^{1,r})}\nonumber\\
\lesssim &||\delta\qnm||_{\infr}T^\frac1p(||\qn||_{\infinf}^2+||\qnm||_{\infinf}^2)\nonumber\\
& +||\nabla\delta\qnm||_{\pinf}(||\qn||_{\infinf}^2+||\qnm||_{\infinf}^2)\nonumber\\
&+||\delta\qnm||_{\infinf}(||\nabla\qn||_{\pinf}+||\nabla\qnm||_{\pinf})\nonumber\\
&\cdot(||\qn||_{\infr}+||\qnm||_{\infr})\nonumber\\
\lesssim &36C^2U_0^2(T^\frac 1p||\delta\qnm||_{\infr}+||\nabla\delta\qnm||_{\pinf}+T^{\frac 12(2-\frac 3r)}||\delta\qnm||_{\infinf})\nonumber\\
\lesssim& (T^\frac 1p+T^{\frac 12(2-\frac 3r)})\delta U^{n-1}(t).\label{3.16}
\end{align}

By Theorem \ref{thm2.1} and Theorem \ref{thm2.2}, and the estimate above, one obtains that
\begin{equation}\label{eqn3.8}
\delta U^n(t)\le C(T^{\frac 12-\frac 3{2q}}+T^{\frac 12(1-\frac 3r)}+T^{\frac 3{2p}(\frac 1q-\frac 1r)}+T^{\frac 12(2-\frac 3r)}+T^\frac 1p)\delta U^{n-1}(t).
\end{equation}
Since the index of $T$ is all positive, one could choose a positive $T_0$ sufficiently small such that
\begin{equation*}
\delta U^n(t)\le\dfrac 12\delta U^{n-1}(t)\quad\forall t\in [0,T_0].
\end{equation*}
$T_0$ only depends on the initial value $U_0$ and $a,b,c,p,q,r$.

It is clear now that the sequence $\{(u^n,\qn,p^n)\}_{n\in\mathbb{N}}$ is a Cauchy sequence in Banach space $M^{p,q,r}_{T_0}$.\\

\noindent
\textit{Step 4. The limit is a strong solution}.\\

Suppose $(u^n,\qn,p^n)\longrightarrow (u,Q,p)$ in $M^{p,q,r}_{T_0}$, obviously the left hand side of (\ref{eqn3.1}) converges to the corresponding ones in (\ref{eqn1.6}) in $L^p(0,T_0,L^q(U))$ and $L^p(0,T_0;W^{1,r}(U))$ respectively. Now we claim that all the nonlinear terms on the right hand side of (\ref{eqn3.1}) also converges to the corresponding ones in (\ref{eqn1.6}) in $L^p(0,T_0;L^q(U))$ and $L^p(0,T_0;W^{1,r}(U))$ respectively. For the brief of the statement, we just take two nonlinear terms with highest order as example. 

For starters, one should notice, as $n\to \infty$
\begin{align*}
\begin{cases}
U^n(t)\longrightarrow U(t)\quad\text{on} [0,T_0],\\
U^n(t)\le 6 CU_0\quad\text{on} [0,T_0].
\end{cases}
\end{align*}
Hence, $U(t)\le 6CU_0$ on $[0,T_0]$.
\begin{align*}
&||\divv(\nabla\qn\odot\nabla\qn-\nabla Q\odot\nabla Q)||_{L^p(0,T_0;L^q)}\\
&\le||\divv(\nabla(\qn-Q)\odot\nabla\qn)||_{L^p(0,T_0;L^q)}+||\divv(\nabla Q\odot\nabla(\qn-Q))||_{L^p(0,T_0;L^q)}\\
&\le||\nabla^2(\qn-Q)||_{L^p(0,T_0;L^\infty)}(||\nabla\qn||_{L^\infty(0,T_0;L^q)}+||\nabla Q||_{L^\infty(0,T_0;L^q)})\\
&\quad+||\nabla(\qn-Q)||_{L^\infty(0,T_0;L^q)}(||\nabla^2\qn||_{L^p(0,T_0;L^\infty)}+||\nabla^2 Q||_{L^p(0,T_0;L^\infty)})\\
&\lesssim ||\qn-Q||_{M^{p,q,r}_{T_0}}T_0^{\frac12(1-\frac 3r)}(U^n(T_0)+U(T_0))\\
&\lesssim 12CU_0||\qn-Q||_{M^{p,q,r}_{T_0}}T_0^{\frac12(1-\frac 3r)}\\
&\longrightarrow 0\quad\text{as} \quad n\to +\infty,
\end{align*}

\begin{align*}
&||\divv(\Delta\qn\qn-\qn\Delta\qn-\Delta Q Q+Q\Delta Q)||_{L^p(0,T_0;L^q)}\\
&\le||\divv(\Delta(\qn-Q)\qn)||_{L^p(0,T_0;L^q)}+||\divv(\Delta Q(\qn-Q))||_{L^p(0,T_0;L^q)}\\
&\quad+||\divv((\qn-Q)\Delta\qn)||_{L^p(0,T_0;L^q)}+||\divv(Q(\Delta(\qn-Q)))||_{L^p(0,T_0;L^q)}\\
&\le||\nabla^3(\qn-Q)||_{L^p(0,T_0;L^q)}(||\qn||_{L^\infty(0,T_0;L^\infty)}+||Q||_{L^\infty(0,T_0;L^\infty)})\\
&\quad+||\Delta(\qn-Q)||_{L^p(0,T_0;L^\infty)}(||\nabla\qn||_{L^\infty(0,T_0;L^q)}+||\nabla Q||_{L^\infty(0,T_0;L^q)})\\
&\quad+||\nabla(\qn-Q)||_{L^\infty(0,T_0;L^q)}(||\Delta Q||_{L^p(0,T_0;L^\infty)}+||\Delta\qn||_{L^p(0,T_0;L^\infty)})\\
&\quad+||\qn-Q||_{L^\infty(0,T_0;L^\infty)}(||\nabla^3 Q||_{L^p(0,T_0;L^q)}+||\nabla^3\qn||_{L^p(0,T_0;L^q)})\\
&\lesssim ||\qn-Q||_{M^{p,q,r}_{T_0}}(U^n(T_0)+U(T_0))(T^{\frac 3{2p}(\frac 1q-\frac 1r)}_0+T^{\frac 12(1-\frac 3r)}_0)\\
&\lesssim 12CU_0||\qn-Q||_{M^{p,q,r}_{T_0}}(T^{\frac 3{2p}(\frac 1q-\frac 1r)}_0+T^{\frac 12(1-\frac 3r)}_0)\\
&\longrightarrow 0\quad\text{as}\quad n\to +\infty.
\end{align*}
where we have used Lemma \ref{lemma2.1} and \ref{lemma2.2}. Similarly one can obtain
\begin{align*}
& u^n\cdot\nabla u^n\longrightarrow u\cdot\nabla u\quad\text{in}\quad L^p(0,T_0;L^q);\\
& u^n\cdot\nabla\qn\longrightarrow u\cdot\nabla Q\quad\text{in}\quad L^p(0,T_0;W^{1,r});\\
& a\qn+b[{\qn}^2-\dfrac{Id}{d}\trace({\qn}^2)]-c\qn\trace({\qn}^2)\\
&\qquad\longrightarrow aQ+b[Q^2-\dfrac{Id}{d}\trace(Q^2)]-cQ\trace(Q^2)\quad\text{in}\quad L^p(0,T_0;W^{1,r})
\end{align*}

\noindent
\textit{Step 5. Uniqueness of the local solution}.\\

The proof of uniqueness is somewhat standard and trivial, we only sketch the proof here.

Let $(u_i,Q_i,p_i)$, $i=1,2$ be two set of solutions to (\ref{eqn1.6}) and (\ref{1.8}) in $M^{p,q}_{T_0}$. Denote their difference by
\begin{equation*}
\delta u:=u_1-u_2,\quad \delta Q:=Q_1-Q_2,\quad \delta p:=p_1-p_2.
\end{equation*}
Then $(\delta u,\delta Q,\delta p)$ satisfies the system
\begin{align*}
\begin{cases}
&\partial_t\delta u-\Delta\delta u+\nabla\delta p=-u_1\cdot\nabla u_1+u_2\cdot\nabla u_2-\divv(\nabla Q_1\odot\nabla Q_1-\nabla Q_2\odot\nabla Q_2)\\
&\quad -\divv(\Delta Q_1Q_1-Q_1\Delta Q_1-\Delta Q_2Q_2+Q_2\Delta Q_2)\\
&\partial_t\delta Q-\Delta\delta Q=-u_1\cdot\nabla Q_1+u_2\cdot\nabla Q_2+\Omega_1Q_1-\Omega_2Q_2-Q_1\Omega_1+Q_2\Omega_2\\
&\quad -a\delta Q+b[Q_1^2-Q_2^2-\dfrac{Id}{d}(\trace(Q_1^2-Q_2^2))]-cQ_1\trace(Q_1^2)+cQ_2\trace(Q_2^2)\\
&\divv \delta u=0,\qquad\int_{\Omega}\delta pdx=0\\
&(\delta u,\delta Q)|_{t=0}=(\delta u,\partial_n\delta Q)|_{\partial U}=(0,0).
\end{cases}
\end{align*}

\begin{equation*}
\begin{aligned}
\delta U(t):=&||\delta u||_{L^\infty(0,t;D^{1-\frac 1p,p}_{A_q})}+||\delta u||_{L^p(0,t;W^{2,q})}+||\partial_t\delta u||_{L^p(0,t;L^q)}\\
&||\nabla\delta p||_{\pq}+||\delta Q||_{L^\infty(0,t;B^{3-\frac 2p}_{r,p})}+||\delta Q||_{\mathcal{W}(0,t)}.
\end{aligned}
\end{equation*}

Repeat the argument of Step 3, one obtains
\begin{equation*}
0\le\delta U(t)\le \dfrac 12 \delta U(t)\quad\text{on}\quad[0,T_0],
\end{equation*}
which implies $\delta U(t)\equiv 0$ on $[0,T_0]$. Thus the uniqueness follows.

\vskip 0.5cm
\section{Global existence of solution with small initial data}

The global existence of strong solution with sufficiently small initial data is proved in this section. Denote $T^*$ the maximal time of existence for $(u,Q,p)$ to (\ref{eqn1.6}) (\ref{1.8}). To this end, one need to control  $T^*$ only in terms of the initial data to extend the local strong solution.

To prepare for the global existence of solution, we prove the low order estimate for the $Q$-tensor first.
\begin{lemma}\label{lemma4.1}
	Suppose $Q\in S^3_{(0)}:=\{Q\in\mathbb{M}^{3\times 3};Q=Q^T;\trace Q=0 \}$, then
	\begin{equation}\label{eqn4.1}
	\trace{(Q^3)}\le \dfrac {3\epsilon}8\trace^2(Q^2)+\dfrac 3{2\epsilon}\trace(Q^2)
	\end{equation}
	for arbitrary $\epsilon>0$.
\end{lemma}
\textit{Proof of the lemma}: 

Since $Q$ is symmetric and trace-free, one may assume eigenvalues of $Q$ are $x$, $y$ and $-x-y$, where $x$, $y\in\mathbb{R}$. Then direct computation shows
\begin{equation*}
\trace(Q^2)=x^2+y^2+(x+y)^2.
\end{equation*}

Hence,
\begin{align*}
\trace(Q^3)=& x^3+y^3-(x+y)^3\\
  =&-3xy(x+y)\\
\le& \dfrac 32 (\epsilon (xy)^2+\dfrac 1\epsilon (x+y)^2)\\
\le&\dfrac{3}{2}(\dfrac \epsilon 4((x+y)^2)^2+\dfrac 1\epsilon(x^2+y^2+(x+y)^2))\\
\le& \dfrac{3\epsilon}{8}\trace^2(Q^2)+\dfrac 3{2\epsilon}\trace(Q^2)
\end{align*}
where we used Young's inequality in above.
\hfill $\Box$

Under the suitable condition of the coefficients in the Landau-de Gennes energy, one can derive exponential decay of the lower order energy.

\begin{proposition}\label{prop4.1}
	Suppose $(u,Q,p)$ satisfy the system (\ref{eqn1.6}) and initial boundary condition (\ref{1.8}) a.e. on $U\times\mathbb{R}_+$, if $Q_0\in L^p(U)$ for some $2\le p<\infty$, then
	\begin{equation*}
	||Q(t,\cdot)||_{L^p}\le e^{-Ct}||Q_0||_{L^p},\qquad \forall t\ge 0
	\end{equation*}
	with some positive  constant $C=C(a,b,c,p)$.
\end{proposition}
\textit{Proof of the proposition}:
 
Multiply the first equation in (\ref{eqn1.6}) by $pQ_{\alpha\beta}\trace^{\frac p2-1}(Q^2)$, sum over the greek indices,
\begin{align}\label{eqn4.2}
(\partial_t+u\cdot\nabla)\trace^{\frac p2}(Q^2)=&p\Delta Q_{\alpha\beta}Q_{\alpha\beta}\trace^{\frac p2-1}(Q^2)-pa\trace^{\frac p2}(Q^2)+pb\trace^{\frac p2-1}(Q^2)\trace(Q^3)-pc\trace^{\frac p2+1}(Q^2)
\end{align}
where we have used fact that $(\Omega_{\alpha\gamma}Q_{\gamma\beta}-Q_{\alpha\gamma}\Omega_{\gamma\beta})Q_{\alpha\beta}=0$ due to the symmetry of $Q$.

Integrate (\ref{eqn4.2}) over the domain $U$ and use the boundary conditions to  integrate by part. Use Lemma \ref{lemma4.1} to control $\trace(Q^3)$.
\begin{align}
\dfrac d{dt}\int_{U}\trace^{\frac p2}(Q^2)dx\le&-p\int_{U}\partial_kQ_{\alpha\beta}\partial_kQ_{\alpha\beta}\trace^{\frac p2-1}(Q^2)dx-p(\dfrac p2-1)\int_{U}|\nabla\trace(Q^2)|^2\trace^{\frac p2-2}(Q^2)dx\nonumber\\
&+p(|b|\dfrac{3\epsilon}8-c)\int_{U}\trace^{\frac p2+1}(Q^2)dx+p(|b|\dfrac 3{2\epsilon}-a)\int_{U}\trace^{\frac p2}(Q^2)dx\nonumber\\
\le& -p(a-\frac{9b^2}{16c})\int_{U}\trace^{\frac p2}(Q^2)dx,\label{4.3}
\end{align}
where $\epsilon=\frac{8c}{3|b|}$. In the last inequality, the first two terms on the right hand side  is negative (since $p\ge 2$). And by (\ref{damping}), 
solve the ordinary inequality (\ref{4.3}), the conclusion follows.

\hfill $\Box$

\textit{Proof of Theorem \ref{existence}}

With the proposition above, next we prove the global solution with sufficiently small initial data.

Define
\begin{align*}
H(t):=& ||u||_{L^\infty(0,t;D^{1-\frac 1p,p}_{A_q})}+||u||_{L^p(0,t;W^{2,q})}+||\nabla p||_{L^p(0,t;L^q)}\\
&+||\partial_t u||_{\pq}+||Q||_{L^\infty(0,t;B^{3-\frac 2p}_{r,p})}+||Q||_{\mathcal{W}(0,t)},
\end{align*}
where
\begin{equation}\label{4.1}
\mathcal{W}(0,t):=W^{1,p}(0,t;W^{1,r})\cap L^p(0,t;W^{3,r});
\end{equation}
and
\begin{equation}\label{4.2}
H(0)=H_0:=||u_0||_{D^{1-\frac 1p,p}_{A_q}}+||Q_0||_{B^{3-\frac 2p}_{r,p}\cap L^r}.
\end{equation}

The by the properties of $L^p$ space, $H(t)$ is non-decreasing and absolutely continuous on $[0,T^*)$. By the maximal regularities of parabolic equations and Stokes operator, as Theorem \ref{thm2.1} and Theorem \ref{thm2.2},
\begin{equation}\label{eqn 4.3}
\begin{aligned}
H(t)\le & C (H_0+\underbrace{||u\cdot\nabla u||_{L^p(0,t;L^q)}}_{\mathcal{K}_1}+\underbrace{||\divv(\nabla Q\odot\nabla Q)||_{L^p(0,t;L^q)}}_{\mathcal{K}_2}\\
&+\underbrace{||\divv(\Delta Q Q-Q\Delta Q)||_{L^p(0,t;L^q)}}_{\mathcal{K}_3}+\underbrace{||u\cdot\nabla Q||_{L^p(0,t;W^{1,r})}}_{\mathcal{K}_4}\\
&\underbrace{||-aQ+b[Q^2-\dfrac{Id}{d}\trace{(Q^2)}]-cQ\trace{(Q^2)}||_{L^p(0,t;W^{1,r})}}_{\mathcal{K}_5}).
\end{aligned}
\end{equation}

By the same method before in Section 4, using Lemma \ref{lemma2.1} and \ref{lemma2.2}, one obtains
\begin{equation}\label{5.7}
\begin{aligned}
H(t)\le & C(T^{\frac 12(1-\frac 3 r)}+T^{\frac 12(1-\frac 3q)}+T^{\frac 3{2p}(\frac 1 q-\frac 1 r)})[H(t)]^2\\
& +C(T^\frac 1pH(t)+T^{\frac 12(2-\frac 3 r)}[H(t)]^3)+H_0\end{aligned}
\end{equation}

Suppose $T$ is the first (smallest) time such that
\begin{equation*}
H(T)= 6CU_0.
\end{equation*} 
Such $T<T^*$ can be found since that $H(t)$ is nondecreasing and continuous in time. Then
\begin{equation*}
H(t)<H(T)=6CH_0,\quad\forall t\in[0,T).
\end{equation*}

Combine this to (\ref{5.7}),
\begin{equation}\label{5.8}
\begin{aligned}
5\le & C(T^{\frac 12(1-\frac 3 r)}+T^{\frac 12(1-\frac 3q)}+T^{\frac 3{2p}(\frac 1 q-\frac 1 r)})H_0\\
& +C(T^\frac 1p+T^{\frac 12(2-\frac 3r)}[H_0]^2).
\end{aligned}
\end{equation}
Observing this inequality, it implies that the maximal existing time $T^*>T$ goes to infinity as the initial data approaches zero. More explicitly, we show next that for sufficiently small initial data, the solution exists globally in time.

\begin{equation}\label{4.6}
\begin{aligned}
\mathcal{K}_1\le& ||u||_{L^\infty(0,t;L^q)}||\nabla u||_{L^p(0,t;L^q)}\\
\le& C ||u||_{L^\infty(0,t;D^{1-\frac 1p,p}_{A_q})}||u||_{L^p(0,t;W^{2,q})}\\
\le& CH^2(t),
\end{aligned}
\end{equation}
\begin{equation}\label{4.7}
\begin{aligned}
\mathcal{K}_2\le& C||\nabla Q||_{L^\infty(0,t;L^q)}||\nabla^2 Q||_{L^p(0,t;L^\infty)}\\
\le& C||Q||_{L^\infty(0,t;W^{1,q})}||Q||_{L^p(0,t;W^{3,r})}\\
\le& C||Q||_{L^\infty(0,t;B^{3-\frac 2p}_{r,p})}||Q||_{\mathcal{W}(0,t)}\\
\le& CH^2(t),
\end{aligned}
\end{equation}
\begin{equation}\label{4.8}
\mathcal{K}_3\le2||\nabla^3QQ||_{L^p(0,t;L^q)}+2||\nabla Q\nabla^2Q||_{L^p(0,t;L^q)}.
\end{equation}

Estimate one by one:
\begin{align*}
||Q\nabla^3Q||_{L^p(0,t;L^q)}\le & C ||Q||_{L^\infty(0,t;L^\infty)}||\nabla^3 Q||_{L^p(0,t;L^q)}\\
\le& C ||Q||_{L^\infty(0,t;B^{3-\frac 2p}_{r,p})}||Q||_{L^p(0,t;W^{3,r})}\\
\le & CH^2(t)
\end{align*}

On the other hand,
\begin{align*}
||\nabla ^2Q\nabla Q||_{L^p(0,t;L^q)}\le &C||\nabla Q||_{L^\infty(0,t;L^q)}||\nabla^2Q||_{L^p(0,t;L^\infty)}\\
\le&C||Q||_{L^\infty(0,t;W^{1,q})}||\nabla^2 Q||_{L^p(0,t;W^{1,q})}\\
\le&C||Q||_{L^\infty(0,t;B^{3-\frac 2p}_{r,p})}||Q||_{L^p(0,t;W^{3,r})}\\
\le&CH^2(t).
\end{align*}
Combine the estimates of the above two terms,
\begin{equation}\label{4.9}
\mathcal{K}_3\le CH^2(t).
\end{equation}

\begin{equation}\label{4.10}
\begin{aligned}
\mathcal{K}_4\le&||u||_{L^\infty(0,t;L^q)}||\nabla Q||_{L^p(0,t;L^\infty)}\\
\le& C ||u||_{L^\infty(0,t;D^{1-\frac 1p,p}_{A_q})}||\nabla Q||_{L^p(0,t;W^{1,r})}\\
\le& CH^2(t).
\end{aligned}
\end{equation}

For the last term, we need to take advantage the decay of lower order norm.
First of all, by Proposition \ref{prop4.1} we have the following
\begin{equation}
||Q||_{L^\infty(0,\infty;L^r)}\le H_0;
\end{equation}
Moreover,
\begin{align}
||Q||_{L^p(0,t;L^r)}= &C(\int_0^t||Q(s,\cdot)||_{L^r}^pds)^{\frac 1p}\nonumber\\
\le&C||Q_0||_{L^r}(\int_0^te^{-pCs}ds)^{\frac 1p}\nonumber\\
\le&C||Q_0||_{L^r}.\nonumber\\
\le & CH_0\label{4.12}
\end{align}
\begin{align}
||\nabla Q||_{\pr}\le &C||Q||_{\pr}^{\frac 12}||\nabla^2 Q||_{\pr}^{\frac 12}\nonumber\\
\le &CH_0^{\frac 12}H^{\frac 12}(t),\label{5.18}
\end{align}

\begin{align}
||\nabla Q||_{L^p(0,t;L^\infty)}\le &C||Q||_{L^\infty(0,t;L^r)}^{\frac 23}||\nabla^3 Q||_{L^p(0,t;L^r)}^{\frac 13}\nonumber\\
\le &CH_0^{\frac 23}H^{\frac 13}(t) \label{5.17}
\end{align}
With the help of above estimate, the last term can be controlled as follows,
\begin{align}
||\mathcal{K}_5||_{L^p(0,t;W^{1,r})}\lesssim& ||Q||_{\pr}+||Q||_{\infr}||Q||_{\pinf}\nonumber\\
&+||Q||_{\infr}||Q||_{\pinf}||Q||_{\infinf}+||\nabla Q||_{\pr}\nonumber\\
&+||\nabla Q||_{\pinf}||Q||_{\infr}+||\nabla Q||_{\pinf}||Q||_{\infr}||Q||_{\infinf}\nonumber\\
\lesssim & H_0+H_0H(t)+H_0H^2(t)+H_0^{\frac 12}H^{\frac 12}(t)+H_0^{\frac 53}H^{\frac 13}(t)+H_0^{\frac 53}H^{\frac 43}(t)\nonumber\\
\lesssim & (1+\epsilon^{-1})H_0+H_0^2+\epsilon^{-1}H_0^{\frac 52}+\epsilon H(t)+(1+H^{5}_0)H^2(t)\label{5.20}
\end{align}
Substitute (\ref{4.6}--\ref{5.20}) into (\ref{4.1}), and choose $C\epsilon<< 1$,
\begin{equation}\label{4.15}
H(t)\le C(H_0+H_0^2+H_0^{\frac 52}+(1+H_0^5)H^2(t)).
\end{equation}
Now take $H_0$ sufficiently small such that
\begin{equation}\label{5.20}
H_0\le \delta_0:=\min\{\dfrac 1{24C^2},1\}.
\end{equation}
Then, directly from (\ref{4.15}), the continuity of $H(t)$ and $H(0)=H_0\le \delta_0$,
\begin{equation}
H(t)\le\dfrac{1-\sqrt{1-24C^2H_0}}{2C(1+H_0^5)}\le\dfrac 1{2C}<\infty,
\end{equation}
for all $t\in[0,T^*)$, which implies $||(u,Q,p)||_{M^{p,q,r}_{T^*}}<\infty$. Thus by the local existence in previous section, one can extend the solution on $[0,T^*)$ to some larger interval $[0,T^*+\delta)$ for some $\delta>0$, which contradicts to the assumption that $T^*$ is the maximal existence time. Hence when the initial conditions satisfy (\ref{5.20}), the local strong solution obtained in previous section is indeed global in time.
	
	The proof of the existence to the global strong solution is finished.
	
\hfill
$\Box$
	
\begin{remark}
(i) The small coefficient $\epsilon$ of $H(t)$ in \eqref{5.20} is crucial in our method, and it's the main reason why we impose condition \eqref{damping}. Otherwise, $H(t)$ cannot be absorbed by the left hand side, then \eqref{4.15} could be trivially true without the uniform bound of $H(t)$.\\

(ii) Notice that the one order energy $H(t)$ solely comes from $aQ$, the linear term. Which is not surprising, in this strong solution framework, nonlinear effects serve as a ``good effect", providing extra decay compared to the linear term. Hence as put in the introduction, if $a=0$ then our result will also follow.
\end{remark}

\end{document}